\documentclass[11pt]{article} \usepackage[dvips]{epsfig}
\usepackage{latexsym}
\usepackage{amsfonts,amsmath,amssymb}
\newtheorem{corollary}{Corollary}
\newtheorem{lemma}{Lemma}
\newtheorem{example}{Example}
\newtheorem{theorem}{Theorem}
\newtheorem{proposition}{Proposition}
\newtheorem{question}{Question}
\newcommand{\qed}{\mbox{$\Diamond$}\vspace{\baselineskip}}
\newenvironment{proof}{\noindent{\bf Proof:}}{\qed}

\begin{document}
\author{Mikl\'os B\'ona\\
        Department of Mathematics\\
University of Florida\\
Gainesville FL 32611-8105\\
USA  }

\title{On a Balanced Property of Compositions}
\maketitle

\begin{abstract} Let $S$ be a finite set of positive integers with 
largest element $m$. Let us randomly
select a composition $a$ of the integer $n$ with parts in $S$, and let $m(a)$
be the multiplicity of $m$ as a part of $a$. Let $0\leq r<q$ be integers, with
$q\geq 2$, 
and let $p_{n,r}$ be the probability that $m(a)$ is congruent to $r$
modulo $q$.  We show that if $S$ satisfies 
a certain simple condition, then $\lim_{n\rightarrow \infty} p_{n,r}
=1/q$. In fact, we show that an obvious necessary condition on $S$ turns out
to be sufficient. 
\end{abstract}

\section{Introduction} A {\em composition} of the positive integer $n$
is a sequence $(a_1,a_2,\cdots ,a_k)$ of positive integers so that
$\sum_{i=1}^k a_i=n$. The $a_i$ are called the {\em parts} of the composition.
It is well-known \cite{bonaw}
 that the number of compositions of $n$ into $k$ parts
is ${n-1\choose k-1}$.
From this fact, it is possible to prove the following. 
Let $0\leq r<q$ be integers, with
$q\geq 2$, and let $P_{n,r}$ be the probability of the event that
the number of parts of a randomly selected composition of $n$ is congruent to
$r$ modulo $q$. Then  $\lim_{n\rightarrow \infty} P_{n,r}
=1/q$. In other words, 
\[\lim_{n\rightarrow \infty}
 \frac{\sum_{i=0}^{\lfloor (n-1)/q \rfloor} {n-1\choose iq+r}}{2^{n-1}}=
\frac{1}{q}.\]
For $q=2$, this follows from the well-known fact that the number of 
even-sized subsets of a non-empty
finite set is equal to the number of odd-sized
subsets of that set. For $q>2$, the statement can be proved, for example,
by the method we will use in this paper. The special cases of $q=3$ and $q=4$
appear as Exercises 4.41 and 4.42 in \cite{bonaw}.
In other words, all residue classes
are equally likely to occur. We will refer to this phenomenon by saying
that the number of part sizes of a randomly selected composition of $n$ is
{\em balanced}.   

Now let us impose a restriction on the {\em part sizes} of the compositions
of $n$ that form our sample space by requiring that all part sizes come
from a finite set $S$. Is it still true that the number of part
 sizes of a randomly selected composition of $n$ is balanced? That will 
certainly depend on the restriction we impose on the part sizes. 
For instance, if the set $S$ of allowed parts consists of odd numbers only,
then the number of part sizes will not be balanced. Indeed, if $q=2$, 
then $P_{n,0}=1$ if $n$ is even, and $P_{n,0}=0$ if $n$ is odd. 
Let $m$ be the largest element of the set $S$ of allowed parts. It turns
out that it is easier to (directly) 
work with the multiplicity $m(a)$ of $m$ as a part of the randomly
selected composition $a$ than with its number of parts. 
For the special case when $S$ has only two elements, the
 results obtained for $m(a)$ can then  be translated back to results on
the number of parts of $a$.

In this paper, we prove that if $S$ satisfies a certain obviously necessary
condition, then the parameter $m(a)$ is balanced, as described in the
 abstract.
That is, the remainder of $m(a)$ modulo $q$ is equally likely to take all
possible values.

\section{The Strategy}

Let $S=\{s_1,s_2,\cdots ,s_k=m\}$ be a finite set of positive integers
with at least two elements. Let us assume without loss of generality
 that no integer larger than 1
divides all $k$ elements of $S$.  Clearly, if $s_1,s_2, \cdots ,s_{k-1}$ are
all divisible by a certain prime $h>1$, then the multiplicity
$m(a)$ of $m=s_k$ as a part of a composition $a$ of $n$ is restricted.
Indeed, $n-m(a)m$ must be divisible by $h$. In particular, if $n$ is divisible
by $h$, then $m(a)m$ is  divisible by $h$, and so $m(a)$ must be divisible
by $h$. Therefore, the parameter $m(a)$ is not balanced. Indeed, if $n$
is divisible by $h$, and we choose
$q=h$, then $p_{n,0}=1$, and $p_{n,r}=0$ for $r\neq 0$, while if $n$ is
not divisible by $h$ and $q=h$, then $p_{n,0}=0$. 

So for $m(a)$ to be a balanced parameter, it is {\em necessary} for $S$ to 
satisfy the condition that its smallest $k-1$ elements do not have 
a proper common divisor. In the rest of this paper, we prove that this
condition is at the same time {\em sufficient} for $m(a)$ to be balanced.

{\em Unless otherwise stated}, let $S$ be a finite set of positive integers
 with at least two elements, and let $S=\{s_1,s_2,\cdots ,s_k=m\}$, where 
the $s_i$ are listed in increasing order. So $m$ is the largest element of
$S$. {\em Unless otherwise stated}, let us also assume  that no integer
larger than 1 is a divisor of all of $s_1,s_2, \cdots ,s_{k-1}$. 
Note that this means that if $|S|=2$, then $s_1=1$. 

 For a fixed positive integer $n$, let
$A_{S,n}(x)$ be the ordinary generating function of all compositions of the
integer $n$ into parts in $S$ according to their number of parts equal to
$m$. In other words, \[A _{S,n}(x)=\sum_a x^{m(a)}=
\sum_{d}a_{S,n,d}x^d,\] where
$a$ ranges over all compositions of $n$ into parts in $S$, and $m(a)$ is
the multiplicity of $m$ as a part in $a$. On the far right, $a_{S,n,d}$ is
the number of compositions $a$ of $n$ into parts in $S$ so that $m(a)=d$.

\begin{example} Let $S=\{1,3\}$. Then the first few polynomials $A_n(x)=
A_{S,n}(x)$
are as follows. 
\begin{itemize} 
\item $A_0(x)=A_1(x)=A_2(x)=1$,
\item $A_3(x)=x$, $A_4(x)=2x+1$, $A_5(x)=3x+1$. 
\item $A_6(x)=x^2+4x+1$, $A_7(x)=3x^2+5x+1$, $A_8(x)=6x^2+6x+1$.
\end{itemize}
\end{example}

Let $q\geq 2$ be a positive integer, and let $0\leq r\leq q-1$.
 Let $A_{S,n,r}$ be the number of compositions $a$ of  $n$ with parts in $S$
so that $m(a)$ is congruent
 to $r$ modulo $q$. So
\[A_{S,n,r}=a_{S,n,r}+a_{S,n,q+r}+\cdots +a_{S,n,\lfloor n/q \rfloor q+r} .\]
 In order to simplify the presentations of our results,
we will first discuss the special case when $n$ is divisible by $q$ and
$r=0$. Let $w$ be a primitive $q$th root of unity. Then
\begin{eqnarray}  \label{unitroots} 
\sum_{t=0}^{q-1} A_{S,n}(w^t) & = & \sum_{t=0}^{q-1} 
\sum_{d= 0}^{n/m} a_{S,n,d} w^{td} \\
 & = & \sum_{d=0}^{n/m} \sum_{t=0}^{q-1} 
a_{S,n,d} w^{td} \\
 & = & \sum_{d=0}^{n/m} a_{S,n,d} \sum_{t=0}^{q-1} w^{td}.
\end{eqnarray}

Using the summation formula of a geometric progression, we get that
\[  \sum_{t=0}^{q-1} (w^{d})^t= \left\{ \begin{array}{l@{\ }l}
0 \hbox{ if $w^d\neq 1$, that is,  if $q\nmid d$}, \\
\\ q \hbox{ if $w^d =1$, that is, if $q|d$.
}
\end{array}\right.
\]
Therefore, (\ref{unitroots}) reduces to 
\[
\sum_{t=0}^{q-1} A_{S,n}(w^t)  = 
q\cdot \sum_{j=1}^{n/q} a_{S,n,jq} 
 =  qA_{S,n,0},\]
\begin{equation}
\label{forma0} \frac{1}{q} \sum_{t=0}^{q-1} A_{S,n}(w^t)=A_{S,n,0}.
 \end{equation}

So in order to find the approximate value of $A_{S,n,0}$, it suffices to
find the approximate values of $A_n(w^t)$, for $0\leq t\leq r-1$, and
for a primitive root of unity $w$. The number $A_{S,n}$
 of {\em all} compositions
of $n$ into parts in $S$ is equal to $A_n(1)$, so we will need that value
as well, in order to compute the ratio $A_{S,n,0}/A_{S,n}$. 

Finally, note that if $n$ is not divisible by $q$, but $r=0$, 
then the same argument 
applies, and $\frac{1}{q} \sum_{t=0}^{q-1} A_{S,n}(w^t)
=A_{S,n,0}$ still holds. 
If $r\neq 0$, then instead of computing $\sum_{t=0}^{q-1} A_{S,n}(w^t)$, we 
compute \[T_n(w)=\sum_{t=0}^{q-1} A_{S,n}(w^t)w^{-rt}
=\sum_{d= 0}^{n/m} a_{S,n,d} \sum_{t=0}^{q-1} w^{t(d-r)}.\]
This shows that the coefficient of $w^k$ in $T_n(w)$ is 0, unless
$w^{d-r}=1$, that is, unless $d-r$ is divisible by $q$. If $d-r$ is divisible
by $q$, then this coefficient is $q$. This shows again that
 \begin{equation}
\label{generalc} \frac{1}{q} \sum_{t=0}^{q-1} A_{S,n}(w^t)w^{-tr}=A_{S,n,r}. 
\end{equation}
Therefore, computing $A_{S,n}(w^t)$
 will be useful in the general case as well.

\section{Linear Recurrence Relations}
In order to compute the values of $A_{S,n}(x)$ for various values of $x$,
we can keep $x$ fixed, and let $n$ grow. 
The connection among the polynomials $A_{S,n}(x)$
 is explained by the following Proposition.

\begin{proposition} Let $S=\{s_1,s_2,\cdots ,s_k=m\}$, with $k\geq 2$. Then
for all positive integers $n\geq 3$, the polynomials $A_{S,n}(x)$ satisfy the
recurrence relation 
\begin{equation} \label{recur}
A_{S,n}(x)=A_{S,n-s_1}(x)+A_{S,n-s_2}(x)+\cdots +A_{S,n-s_{k-1}}(x)+
xA_{S,n-m}(x).
\end{equation}
\end{proposition}

\begin{proof}
Let $a$ be a composition of $n$ with parts in $S$. If the first part of
$a$ is $s_i$, for some $i\in [1,k-1]$, then the rest of $a$ forms a 
composition of $n-s_i$ with parts in $S$ in which the multiplicity of $m$ as
a part is still $m(a)$. These compositions of $n$ 
 are counted by $A_{S,n-s_i}(x)$.
If the first part of $a$ is $m$, then the rest of $a$ forms a 
composition of $n-m$ with parts in $S$ in which the multiplicity of $m$ as
a part is $m(a)-1$. These compositions of $n$  are counted by 
$xA_{S,n-m}(x)$.
\end{proof}

\begin{example} If $S=\{1,3\}$, then (\ref{recur}) reduces to 
\begin{equation} \label{srecur} A_{S,n}(x)= A_{S,n-1}(x)+x A_{S,n-3}(x).
\end{equation}
\end{example}

For a {\em fixed} real number $x$, the recurrence relation (\ref{recur})
becomes a recurrence relation on real numbers. The solutions of such
 recurrence relations are described by the following well-known theorem.
(See, for instance, \cite{rosen}, Section 7.2.
)
\begin{theorem} \label{recthe}
Let \begin{equation}
\label{grec} a_n=c_1a_{n-1}+c_2a_{n-2}+\cdots +c_ka_{n-k} \end{equation}
 be a recurrence
relation,
where the $c_i$ are complex constants. Let $\alpha_1,\alpha_2,\cdots,
\alpha_t$ be the distinct roots of the characteristic
equation \begin{equation} \label{chareq}
z^k-c_1z^{k-1}-c_2z^{k-2}-\cdots -c_k=0,\end{equation}
and let $M_i$ be the multiplicity of $\alpha_i$.
Then the sequence $a_0,a_1,\cdots $ of complex numbers satisfies (\ref{grec})
if and only if there exist constants $b_1,b_2, \cdots ,b_k$ so that for 
all $n\geq 0$, we have
\begin{equation}
\label{fgrec} a_n=b_1\alpha_1^n+b_2 n\alpha_1^n+ \cdots 
+b_{M_1}n^{M_1-1}\alpha_1^n+ b_{M_1+1}\alpha_2^n, \cdots
\end{equation}  \begin{equation}
\cdots +b_{M_1+M_2}n^{M_2-1}\alpha_2^n, \cdots, \cdots 
+b_kn^{M_k-1}\alpha_k^n.\end{equation}

In other words, the solutions of (\ref{grec}) form a $k$-dimensional vector
space.
\end{theorem}

We will need the following consequence of Theorem \ref{recthe}.
\begin{corollary} \label{nonzero}
Let us assume that the sequence $\{a_n\}$ is a solution of (\ref{grec}) and
that there is no linear recurrence relation of a degree less than
$k$ that is satisfied by $ \{a_n\}$. Let us further assume that the
characteristic equation (\ref{chareq}) of (\ref{grec}) has a unique
root $\alpha_1$ of largest modulus. Then there is a {\em nonzero}
 constant $C$ so that
\[a_n=C\alpha_1^n + o(\alpha_1^n).\] 
\end{corollary}

\begin{proof}
As $\{a_n\}$ does not satisfy a recurrence relation of a degree less than
$k$, we must have $c_1\neq 0$. As $|\alpha_1|>|\alpha_i|$ for $i\neq 1$, 
the statement follows. 
\end{proof}

Let us now apply Theorem \ref{recthe} to find 
 the solution of (\ref{recur}) for a fixed $x$. The characteristic
equation of (\ref{recur}) is
 
\begin{equation} \label{polyeq} 
f(z)=z^{m}-\sum_{i=1}^{k-1}  z^{m-s_{i}} - x
=0.\end{equation}

As explained in Section 2, we will need to compute $A_{S,n}(1)$ and also,
$A_{S,n}(w^t)$ for the case when $w\neq 1$ is a $q$th primitive
root of unity.  To that
end, we need to find the roots of the corresponding characteristic
equations. That is, we will compare the root of largest modulus
 of the characteristic 
equation 
\begin{equation} \label{forone}
f_1(z)=
z^{m}-\sum_{i=1}^{k-1}  z^{m-s_{i}} - 1=0
\end{equation}
and the root(s) of the largest modulus of the characteristic equation
\begin{equation} \label{forw}
f_w(z)=
z^{m}-\sum_{i=1}^{k-1}  z^{m-s_{i}} - w=0
\end{equation}

The following lemma, helping to compute root of largest
modulus of $f_1(z)$, 
 is a special case of Exercise III.16 in \cite{polya}.

\begin{lemma} \label{first}
 The polynomial  $f_1(z)=z^{m}-\sum_{i=1}^{k-1} 
 z^{m-s_{i}} - 1$
 has a unique positive real root $\alpha$.
\end{lemma}

\begin{proof} Let $\alpha$ be the smallest positive real root of $f_1(z)$.
 We know such a root exists since $f(0)=-1$ and $\lim_{z\rightarrow
\infty}f(z)=\infty$. We claim that then $f$ is strictly monotone increasing
on $[\alpha, \infty )$, implying that $f$ cannot have another positive
real root. Indeed, if $r>1$, then 
\begin{eqnarray*} f(r\alpha)+1 & = & (r\alpha)^{m} - \sum_{i=1}^{k-1} 
 (r\alpha)^{m-s_{i}}\\
& > & r^{m} \left(\alpha^{m}-\sum_{i=1}^{k-1}  \alpha^{m-s_{i}} \right ) \\
& = & r^m \\
& > & 1, 
\end{eqnarray*}
and so $f(r\alpha)>0$.
\end{proof}

Now we address the problem of finding the roots of the characteristic
equation (\ref{polyeq}) in the case when $w\neq 1$ is a root of unity.
It turns out that it suffices to assume that $|w|=1$. (The following
Lemma is similar to Exercise III.17 in \cite{polya}.)

\begin{lemma} Let $\alpha$ be defined as in Lemma \ref{first}. 
Let $w$ be any complex number satisfying  $w\neq 1$ and $|w|=1$.
Then all  roots of the polynomial $f_w(z)$ are of smaller modulus than
 $\alpha$.
\end{lemma}

\begin{proof} Let $y$ be a root of $f$. Then 
\begin{eqnarray*} |y|^m & = & \left |w + \sum_{i=1}^{k-1} y^{m-s_i} 
\right |\\
& \leq & 1+  \sum_{i=1}^{k-1} \left | y^{m-s_i} \right |.
\end{eqnarray*}

 Therefore, $f(|y|)\leq 0$. This implies that $|y|\leq \alpha$ since we have
seen in the proof of Lemma \ref{first} that $f(t)>0$ if $t>\alpha$.

Furthermore, in the last displayed inequality, the inequality is strict
unless for all $i$ so that $1\leq i\leq k-1$,
 the complex numbers $y^{m-s_i}$ have the same argument as $w$, and that
argument is the same as the argument of $y^m$. That happens only if
the complex numbers $y^{s_1},y^{s_2-s_1},\cdots y_{s_{k-1}-s_{k-2}}$ all
have argument 0, that is, when these numbers are positive real numbers.
However, that happens precisely when $s_1, s_2-s_1, \cdots, s_{k-1}-s_{k-2}$
are all multiples of the multiplicative order $o_y$
 of $y/|y|$ as a complex number.
That implies that $s_1,s_2,\cdots ,s_{k-1}$ are all divisible by $o_y$, 
contradicting our hypothesis on $S$. 
\end{proof}

The previous two lemmas show that the largest root of the
 characteristic equation for the
sequence $\{A_{S,n}(1)\}_{n\geq 0}$ is larger than the largest root(s) of
the characteristic equation for the sequence $\{A_{S,n}(w)\}_{n\geq 0}$ for 
any complex number $w\neq 1$ with absolute value 1. Given formula 
(\ref{fgrec}), in order to see that
the first sequence indeed grows faster than the second, 
all we need to show is that the {\em coefficient} $b_1$ of $\alpha^n$
in (\ref{grec}) is not 0. (Here $\alpha$, the largest root of $f_1(z)$, 
plays the role of $\alpha_1$ in (\ref{fgrec})). This is the content of the
next lemma.

\begin{lemma} \label{noshorter} Let $S=\{s_1,s_2,\cdots ,s_{k-1}\}$
 be any finite set of positive
integers (so for this Lemma, we do not require that $s_1,s_2,\cdots s_{k-1}$
 do not have a proper common divisor).
Then the sequence $\{A_{S,n}\}_{n\geq 0}=\{A_{S,n}(1)\}_{n\geq 0}$
 does not satisfies a linear recurrence
relation with constant coefficients and
less than $|S|+1$ terms. In other words, if $|S|=k$, then
 there do not exist
constants $c_2,c_3,\cdots ,c_k$ and positive integers $j_1,j_2,\cdots ,
j_{k-1}$ so that for all $n\geq 0$, 
\[A_{S,n}=\sum_{i=1}^{k-1} c_i A_{S,n-j_i}.\]
\end{lemma}

\begin{proof}
Let us assume that $S$ is a minimal counterexample. It is then straightforward
to verify that $|S|>2$. Let $S'=S-m$, that is, the set obtained from 
$S$ by removing the largest element of $S$. Then
\begin{equation} \label{reduction} 
A_{S',n}=\sum_{i=1}^{k-1} c_i' A_{S',n-s_i},\end{equation}
and there is no shorter recurrence satisfied by $\{A_{S',n}\}$.

Now crucially, $A_{S',n}=A_{S,n}$ for all $n$ satisfying $0\leq n< m$.
So these sequences agree in $m-1\geq k-1$ values. So if $\{A_{S,n}\}$
satisfied a linear recurrence relation of degree $k-1$, that would have
to be the recurrence relation (\ref{reduction}). Indeed, by Theorem 
\ref{recthe}, the solutions of (\ref{reduction}) form a $k$-dimensional 
vector space, so knowing $k$ elements of a solution determines the whole
solution.  
However, $\{A_{S,n}\}$ does not satisfy (\ref{reduction}) since 
$A_{S,m}=A_{S',m}+1\neq A_{S',m}$, where the difference is caused by the
one-part composition $m$. 
\end{proof}

Now we are in position to express the growth rate of $A_{S,n}=A_{S,n}(1)$. 
\begin{proposition}
Let $\alpha$ be defined as in Lemma \ref{first}. Then
\[A_{S,n}(1)=C \alpha^n + o(\alpha^n),\] for some  nonzero
constant $C$. 
\end{proposition}

\begin{proof} Immediate from Corollary \ref{nonzero} and 
Lemma \ref{noshorter}.
\end{proof}

We can now compare the growth rates of $A_{S,n}(w)$ and $A_{S,n}(1)$.

\begin{lemma} \label{grrate}
Let $w\neq 1$ be any complex number so that $|w|=1$. Then
\[\lim_{n\rightarrow \infty} \frac{A_{S,n}(w)}{A_{S,n}(1)} =0. \]
\end{lemma}

\begin{proof}
Lemma (\ref{first}) shows that the unique positive root of the characteristic
equation (\ref{forone}) is larger than the absolute value of all
roots of the characteristic
equation (\ref{forw}). Therefore, $A_{S,n}(w)=O(n^k\beta^k)$, with
$\beta < \alpha$. 
\end{proof}

Finally, we can use the results of this section to prove the balanced 
properties of the numbers $A_{S,n,r}$.

\begin{theorem} 
Let $q\geq 2$ be an integer, and let $r$ be an integer satisfying $0\leq r
\leq q-1$. Then 
\[\lim_{n\rightarrow \infty} p_{n,r}=\lim_{n\rightarrow \infty} 
\frac{A_{S,n,r}}{A_{S,n}}=\frac{1}{q}.\]
\end{theorem}

\begin{proof} Let us first address the case of $r=0$. 
Dividing (\ref{forma0}) by $A_{S,n}$, we get  that 
 \[A_{S,n,0}=\frac{1}{q}
 \sum_{t=0}^{q-1}\frac{A_{S,n}(w^t)}{A_{S,n}}.\] 
However, Lemma \ref{grrate} shows that all but one of the $q$ summands
on the right-hand side converge to 0, and the 
remaining one (the first summand) is equal to 1. 

For general $r$, the only change is that instead of dividing both sides
of (\ref{forma0}) by $A_{S,n}$, we divide both sides of (\ref{generalc})
by $A_{S,n}$. As $|w|=1$, the result follows in the same way. 
\end{proof}

\section{Further Directions}
Let $S=\{1,3\}$. Numerical evidence suggest that for all $n$, the polynomials 
$A_{S,n}(x)$ have real roots only. Furthermore, numerical evidence also
suggests that the sequences of polynomials $\{A_{S,3n+r}(x)\}_{n\geq 0}$ form
a Sturm sequence for each of  $r=0,1,2$. (See \cite{wilfb} for the definition
and importance of Sturm sequences.) This raises the following intriguing
questions. 

\begin{question} For which sets $S$ is it true that  the polynomials 
$A_{S,n}(x)$ have real roots only?
\end{question}

\begin{question} For which sets $S$ is it true that  the polynomials 
$A_{S,n}(x)$ can be partitioned into a few Sturm sequences?
\end{question}

Herb Wilf \cite{wilf} 
proved that the set $S=\{1,2\}$ does have both of these properties. 

If $A_{S,n}(x)$ has real roots only, then its coefficients form a log-concave
(and therefore, unimodal) sequence. (See Chapter 8 of 
\cite{bonaint} for an
introduction into into unimodal and log-concave sequences.)  This raises the
 following questions.

\begin{question} Let us assume that $A_{S,n}(x)$ has real roots only.
Is there a combinatorial proof for the log-concavity of its coefficients?
\end{question}

\begin{question}  Let us assume that $A_{S,n}(x)$ has real roots only.
Where is the peak (or peaks) of the unimodal sequence of its coefficients?
\end{question}

Another interesting question is the following.
\begin{question}
For what $S$ and $n$ does the equality $A_{S,n}(-1)=0$ hold? When it does,
the number of compositions of $n$ with parts in $S$ and with $m(a)$ even
equals the number of compositions of $n$ with parts in $S$ and with $m(a)$
odd. Is there a combinatorial proof of that fact?
\end{question}
 
Finally, our methods rested on the finiteness of $S$, but we can still
ask what can be said for {\em infinite} sets of allowed parts. 

\begin{center} {\bf Acknowledgment} 

I am grateful to Herb Wilf for valuable discussions and advice.
\end{center}


\begin{thebibliography}{99}
\bibitem{bonaw} M. B\'ona, A Walk Through Combinatorics, 2nd edition,
 {\em World Scientific}, 2006.
\bibitem{bonaint}  M. B\'ona, Introduction to Enumerative Combinatorics,
McGraw-Hill, 2007. 
\bibitem{polya} G. P\'olya, G. Szeg\H{o}, {\em Problems and Theorems 
in Analysis}, Springer, 1972.
\bibitem{rosen} K. Rosen, Discrete Mathematics and Its Applications, 
McGraw-Hill, Sixth Edition, 2007. 
\bibitem{wilf} H. Wilf, {\em Personal Communication}, 2006. 
\bibitem{wilfb} H. Wilf, {\em Mathematics for the Physical Sciences},
third edition, Dover, 2006.
\end{thebibliography}
\end{document}